\documentclass[12pt]{amsart}

\usepackage[T1]{fontenc}
\usepackage[utf8]{inputenc}
\usepackage{lmodern}

\usepackage[american]{babel}
\usepackage[babel, final]{microtype}
\usepackage{amssymb}
\usepackage[all]{xy}
\usepackage{amsmath}

\usepackage{ifthen}
\usepackage{ifpdf}
\usepackage{multirow}
\usepackage[shortlabels]{enumitem}
\makeatletter
\newcommand{\mylabel}[2]{#2\def\@currentlabel{#2}\label{#1}}
\makeatother

\ifpdf
  \usepackage[pdftex, dvipsnames]{xcolor}
  \usepackage[pdftex, final]{graphicx}
  \usepackage{caption}
  \usepackage{subcaption}
  \usepackage{pinlabel}
  \usepackage[pdftex,%
    letterpaper,%
    includehead,%
    includefoot,%
    nomarginpar,%
    lmargin=1in,%
    rmargin=1in,%
    tmargin=1in,%
    bmargin=1in,%
  ]{geometry}
  \usepackage[pdftex,%
    final,%
    colorlinks=true,%
    linkcolor=NavyBlue,%
    citecolor=NavyBlue,%
    filecolor=NavyBlue,%
    menucolor=NavyBlue,%
    urlcolor=NavyBlue,%
    bookmarks=true,%
    bookmarksdepth=3,%
    bookmarksnumbered=true,%
    bookmarksopen=true,%
    bookmarksopenlevel=2,%
  ]{hyperref}
  \hypersetup{
    pdftitle={Khovanov homology and Lagrangian cobordisms},
    pdfauthor={Gage Martin, Ina Petkova, Zachary Winkeler},
    pdfsubject={Lagrangian cobordisms},
    pdfkeywords={Lagrangian cobordism, Legendrian knots}
  }
\else
  \usepackage[dvips, dvipsnames]{xcolor}
  \usepackage[dvips, final]{graphicx}
  \usepackage{caption}
  \usepackage{subcaption}
  \usepackage{pinlabel}
  \usepackage[dvips,%
    letterpaper,%
    includehead,%
    includefoot,%
    nomarginpar,%
    lmargin=1in,%
    rmargin=1in,%
    tmargin=1in,%
    bmargin=1in,%
  ]{geometry}
  \usepackage[ps2pdf,%
    final,%
    colorlinks=true,%
    linkcolor=NavyBlue,%
    citecolor=NavyBlue,%
    filecolor=NavyBlue,%
    menucolor=NavyBlue,%
    urlcolor=NavyBlue,%
    bookmarks=true,%
    bookmarksdepth=3,%
    bookmarksnumbered=true,%
    bookmarksopen=true,%
    bookmarksopenlevel=2,%
  ]{hyperref}
  \hypersetup{
    pdftitle={Khovanov homology and Lagrangian cobordisms},
    pdfauthor=Gage Martin, Ina Petkova, Zachary Winkeler},
    pdfsubject={Lagrangian cobordisms},
    pdfkeywords={Lagrangian cobordism, Legendrian knots}
  
\fi

\usepackage{mathdots}
\usepackage{float}
\usepackage[mathscr]{euscript}
\usepackage[normalem]{ulem}
\usepackage{tikz}
\usepackage{tikz-cd}
\usepackage{mathtools}
\usepackage{stmaryrd}
\SetSymbolFont{stmry}{bold}{U}{stmry}{m}{n}
\usetikzlibrary{arrows,automata,fit}
\usetikzlibrary{matrix,arrows}
\usetikzlibrary{decorations.pathreplacing,calligraphy}

\input{macros}

\newcommand{\cross}{\times}                  

\newcommand{\numset}[1]{\mathbb{#1}}

\newcommand{\Z}{\numset{Z}}

\newcommand{\Q}{\numset{Q}}
\newcommand{\R}{\numset{R}}

\newcommand{\xistd}{\xi_{\mathrm{std}}}
\DeclareMathOperator{\tb}{tb}
\DeclareMathOperator{\rot}{r}
\DeclareMathOperator{\self}{sl}

\newcommand*{\alphastd}{\alpha_{\mathrm{std}}}

\newcommand{\std}{\mathrm{std}}

\renewcommand*{\tb}{\mathit{tb}}
\renewcommand*{\rot}{\mathit{r}}

\DeclareMathOperator{\CKh}{CKh}
\DeclareMathOperator{\Kh}{Kh}

\newcommand{\Lee}{\Kh_{Lee}}
\newcommand{\CLee}{\CKh_{Lee}}
\DeclareMathOperator{\gr}{gr}
\newcommand{\complex}{\mathcal{C}}

\newcommand{\Leg}{\Lambda}
\newcommand*{\Legm}{\Leg_-}
\newcommand*{\Legp}{\Leg_+}

\newcommand*{\LKh}{\mathcal L_{\Leg}}
\newcommand*{\LdecKh}{\mathcal L_{\Leg}'}

\newcommand*{\psigen}{\widetilde{\psi}}

\title{Khovanov homology and Lagrangian cobordisms}

\author{Gage Martin}
\address{Department of Mathematics \\ Harvard University \\ Cambridge, MA 02138}
\email{\href{mailto:gagemartin@math.harvard.edu}{gagemartin@math.harvard.edu}}
\urladdr{\url{https://sites.google.com/view/gagemartin/home}}

\author[I. Petkova]{Ina Petkova}
\address{Department of Mathematics \\ Dartmouth College \\ Hanover, NH 03755}
\email{\href{mailto:ina.petkova@dartmouth.edu}{ina.petkova@dartmouth.edu}}
\urladdr{\url{https://math.dartmouth.edu/~ina/}}

\author[Z. Winkeler]{Zachary Winkeler}
\address {Department of Mathematics \\ Colby College \\ Waterville, ME 04901}
\email{\href{mailto:zwinkele@colby.edu}{zwinkele@colby.edu}}
\urladdr{\url{http://zach-winkeler.github.io}}

\begin{document}

\begin{abstract}
    We provide a partial answer to a question of Ekholm, Honda, and K\'{a}lm\'{a}n about the relationship between Khovanov homology and decomposable Lagrangian cobordisms. We also utilize previously defined filtered invariants to give obstructions to decomposable Lagrangian cobordisms from Khovanov homology.
\end{abstract}

\maketitle

\section{Introduction}

Some of the central questions in low-dimensional contact topology concern the classification of Legendrian, as well as transverse, links, and the study of various types of cobordisms between such links. In the simplest case, one studies Legendrian and transverse links in $\mathbb R^3$ with the standard contact structure
  \[
  \xi_{\std}=\ker(\alpha_{\std}),\quad \alpha_{\std}=dz-y\,dx, 
\] 
 and Lagrangian cobordisms (between Legendrian links) in the symplectization 
 \[(\R_t \cross \R^3, d (e^t \alphastd)).\] 

In \cite{olga-psi}, Plamenevskaya introduced an invariant of transverse links coming from Khovanov homology. Given a transverse link $K$, this invariant $\psi(K)$ is an element in the Khovanov homology $\Kh(K)$. Soon after, Ng observed that Khovanov homology can be used to study Legendrian links too, namely, to obtain an upper bound for the maximal Thurston-Bennequin number of a link \cite{NgLenhard2005ALTb}. This bound
\[C = \min \{ k  \mid \bigoplus_{i-j=k}\Kh^{i,j}(K)\neq 0\}\]
can be visually described via the corresponding line $i-j= C$ in the $ji$-plane, which we will refer to as ``Ng's line'' \footnote{In his paper Ng used $i$ to denote the quantum grading and $j$ to denote the homological grading, which is the reverse of contemporary conventions; we follow Ng's convention in this paper.}. Broader connections of Khovanov homology to Legendrian knot theory remain elusive, and past results about Lagrangian cobordisms have had to rely on other tools, such as  tools from symplectic field theory \cite{Cha15:NotSym, ST13, CDGG15, 
  CNS16, Pan17} or Floer homology \cite{BS18, BS21, GJ19, BLW22, JPSWW22}.

Given two Legendrian links $\Legm, \Legp$ in $(\R^3, \xistd)$, Ekholm, Honda, and K\'{a}lm\'{a}n observed that an exact Lagrangian cobordism $L$ from $\Legm$ to $\Legp$ \footnote{Note that Ekholm, Honda, and K\'{a}lm\'{a}n call this ``an exact Lagrangian cobordism $L$ from $\Legp$ to $\Legm$''.} induces a map (defined up to $\pm 1$)
\[F_{mL}\colon \Kh(m (\Legm)) \to \Kh(m(\Legp)),\]
simply by considering the (mirror of the) same cobordism in the smooth category and using functoriality of Khovanov homology. A natural question they then asked is the following: 

\begin{question}[{cf \cite[Question~9.7]{EHK16}}] \label{q:EHK}
Given a Legendrian link $\Leg$, what is the relationship among the set of elements $\LKh \coloneqq \{F_{m L}(\pm 1)\}\subset \Kh(m(\Leg))$ obtained by considering all exact Lagrangian fillings $L$ of $\Leg$, Plamenevskaya's invariant, and Ng's line? Similarly, what about the set $\LdecKh$ where one considers only decomposable Lagrangian fillings? 
\end{question}

At the time of Ekholm, Honda, and K\'{a}lm\'{a}n's paper, given a link bounding a surface, considering surface maps from the ground ring to the Khovanov homology of the link seemed natural; see also \cite{SundbergIsaac2022RKc}. However, more recent developments have shown the power of the dual perspective -- considering surface maps from the Khovanov homology of a link to the ground ring~\cite{hayden_khovanov_2024}. From this dual perspective, we would instead be asking what elements of $\Kh(\Leg)$ are sent to $\pm 1$ by the map $F_L$ for a Lagrangian filling $L$. The maps in both perspectives are simply the smooth cobordism maps established in \cite{J04}; the difference is whether one treats a surface in $\R^3\times \left[0,1\right]$ as a bottom-to-top (as in \cite{EHK16} and \cite{SundbergIsaac2022RKc}) or a top-to-bottom (as in \cite{hayden_khovanov_2024}) cobordism.

In this article, we provide a partial answer to \fullref{q:EHK} in the case of decomposable Lagrangian fillings.

Our strategy involves passing to the category of transverse links and showing that decomposable Lagrangian cobordisms can be deformed to a certain type of ``transverse'' cobordisms between transverse links; this creates the opportunity to utilize a large amount of previous work on Khovanov homology. This previous work lets us completely understand the cobordism maps on Khovanov homology induced by decomposable Lagrangian cobordisms:

\begin{theorem}
\label{thm:KhLagrangian}
    Let $\Leg_-$ and $\Leg_+$ be two Legendrian links with positive transverse pushoffs $K_-$ and $K_+$, and suppose there is a decomposable Lagrangian cobordism $L$ from $\Leg_-$ to $\Leg_+$. Then the induced cobordism map $F_L:  \Kh(\Leg_+)\to \Kh(\Leg_-)$ sends $\psi(K_+)$ to $\pm\psi(K_-)$.
\end{theorem}

Specializing to the case where $\Leg_-$ is empty, and identifying $\Kh(\emptyset) \cong \Z$, we obtain the following: 

\begin{corollary}\label{cor:EHK-answer} 
Let $\Leg\subset (\R^3, \xistd)$ be a Legendrian link with a decomposable Lagrangian filling $L$, and Let $K$ be the positive transverse pushoff for $\Leg$. We have
\[ F_L(\psi(K)) = \pm 1.\]
\end{corollary}

Thinking  of fillings  as bottom-to-top cobordisms instead, we obtain the following:

\begin{theorem}\label{thm:line}
    If $\Leg$ admits a decomposable Lagrangian filling, then the set $\LKh$ is supported in bidegree $(-\tb(\Leg),0)$, and $(-\tb(\Leg),0)$ lies on Ng's line for $\Leg$.
\end{theorem}

As an immediate consequence, we obtain a purely topological obstruction to decomposable Lagrangian fillings:
\begin{corollary}\label{cor:top-obstr}
    If $\Leg$ is a Legendrian knot of topological knot type $K$ that admits a decomposable Lagrangian filling, then 
    \[\overline{\tb}(K) = \min \{ k  \mid \bigoplus_{i-j=k}\Kh^{i,j}(K)\neq 0\}.\]
\end{corollary}

The behavior of the cobordism map asserted in \fullref{thm:KhLagrangian} immediately implies that Plamenevskaya's invariant obstructs the existence of decomposable Lagrangian cobordisms:

\begin{corollary}
\label{cor:Vanish}
    Let $\Leg_-$ and $\Leg_+$ be two Legendrian links with positive transverse pushoffs $K_-$ and $K_+$. If $\psi(K_-) \not=0$ and $\psi(K_+) = 0$, then there is no decomposable Lagrangian cobordism from $\Leg_-$ to $\Leg_+$.
\end{corollary}

Given the relative computability of Khovanov homology, \fullref{cor:Vanish} might initially spark hope in a fast tool for obstructing decomposable Lagrangian cobordisms. Unfortunately, this obstruction can only provide useful information when $\psi$ is not determined solely by smooth invariants and the classical transverse invariant, the self-linking number: 

\begin{theorem}
\label{thm:effectiveness}
Let $\Leg_-$ and $\Leg_+$ be two Legendrian links, and suppose that the positive transverse pushoff $K_-$ of $\Leg_-$ is in the collection of transverse knots $K$ where the conditions ``$\psi(K) \not= 0$'' and ``$s(K) = \self(K) + 1$'' are both true or both false. Then \fullref{cor:Vanish} will not provide a stronger obstruction to the existence of a decomposable Lagrangian cobordism from $\Leg_-$ to $\Leg_+$ than the combination of the classical Legendrian invariants and the $s$-invariant.

\end{theorem}

The question of whether there is a transverse knot $K$ with $s(K) > \self(K)+1$ and $\psi(K)\not=0$ has been a long-standing question at the intersection of Khovanov homology and contact topology. With the previous theorem in mind, it is only possible to get an interesting obstruction from \fullref{cor:Vanish} using such a $K$.

Lipshitz, Ng, and Sarkar~\cite{LipshitzRobert2015Otif} introduced a generalization of the Plamenevskaya invariant $\psi(K)$ using filtrations on the Khovanov homology of $K$. Inspired by results from knot Floer homology~\cite{JPSWW22} and applying similar algebraic techniques we show: 
\begin{theorem}\label{thm:filt}
    The filtered transverse invariants in Khovanov homology are functorial under decomposable Lagrangian cobordisms and provide an obstruction to the existence of such cobordisms.
\end{theorem}
A more precise statement of this result appears in \fullref{thm:Spectral} and \fullref{cor:SpectralVanish}.

\subsection*{Acknowledgments} We thank John Baldwin and Kyle Hayden for helpful conversations. We also thank the anyonymous referee for helpful comments. GM was partially supported by NSF grant DMS-2202841. IP was partially supported by NSF CAREER Grant DMS-2145090.

\section{Decomposable Lagrangian cobordisms and Plamanevskaya's invariant}

The crucial observation needed to understand the interaction between Khovanov homology and decomposable Lagrangian cobordisms is the following:

\begin{proposition}\label{thm:LagrangianToTransverse}
    Let $\Leg_-$ and $\Leg_+$ be two Legendrian links with positive transverse pushoffs $K_-$ and $K_+$, and suppose there is a decomposable Lagrangian cobordism $L$ from $\Leg_-$ to $\Leg_+$. Then there is an ascending cobordism with positive critical points $T$ from $K_-$ to $K_+$ such that $L$ and $T$ are isotopic as smooth cobordisms.
\end{proposition}

A non-precise statement of this observation is that any decomposable Lagrangian cobordism between Legendrian links can be deformed to a ``transverse'' cobordism between transverse links. 

\begin{proof}
We may reduce to the case of an elementary Lagrangian cobordism, as the theorem follows for a general decomposable Lagrangian cobordism simply by breaking this cobordism into elementary pieces.

If the elementary cobordism $L$ from $\Leg_-$ to $\Leg_+$ is a Legendrian isotopy, then the positive transverse pushoffs $K_-$ and $K_+$ are transversely isotopic, and this transverse isotopy induces an ascending cobordism with positive critical points (in fact with no critical points).

If $L$  is a birth (i.e.\ a minimum cobordism), then there is a corresponding cobordism from $K_-$ to $K_+$, obtained by simply deforming $L$ near the boundary. This cobordism is an ascending cobordism with a positive critical point.

Finally, if $L$ is a pinch (i.e.\ a saddle), then consider the pushoffs $K_-$ and $K_+$ in a neighborhood of the saddle. We can construct the cobordism from $K_-$ to $K_+$ by attaching a band with a right-handed half-twist to $K_-$ in that neighborhood; see \fullref{Fig:Band}. This cobordism is an ascending cobordism with a positive critical point.\qedhere
\begin{figure}
  \centering
    \includegraphics[width=0.8\textwidth]{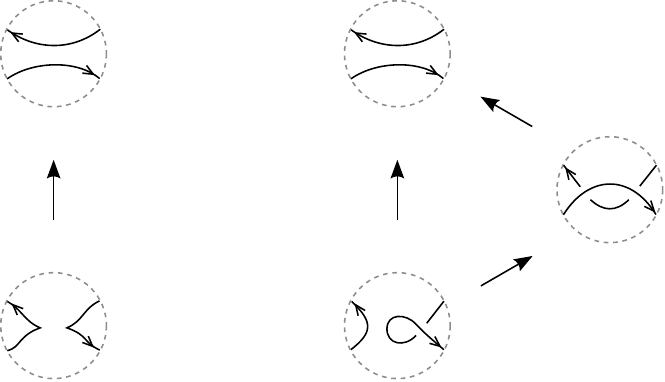}
      \caption{Left: An elementary Lagrangian cobordism for a 1-handle attachment (a pinch) from $\Legm$ to $\Legp$. Right: The corresponding ascending cobordism (a band attachment followed by an isotopy) between the positive transverse pushoffs of $\Legm$ and $\Legp$.}
      \label{Fig:Band}
\end{figure} 
\end{proof}

Having moved from the setting of Legendrians and decomposable Lagrangians to a transverse setting, we are now ready to prove \fullref{thm:KhLagrangian}.

\begin{proof}[Proof of \fullref{thm:KhLagrangian}]
Hayden--Sundberg showed that Plamenevskaya's invariant is functorial under ascending cobordisms with positive critical points \cite[Theorem 1.3]{hayden_khovanov_2024}. 
Combining \fullref{thm:LagrangianToTransverse} with  \cite[Theorem 1.3]{hayden_khovanov_2024} completes the proof.
\end{proof}

\fullref{cor:EHK-answer} follows immediately, by observing that the top-to-bottom cobordism map for a birth from the crossingless unknot to the empty set sends Plamenevskaya's cycle to $\pm 1$. 

\fullref{cor:Vanish}  follows immediately, too. In other words, $\psi$ obstructs the existence of  decomposable Lagrangian cobordisms. 

We next discuss the effectiveness of this obstruction. Specifically, we discuss how the effectiveness of $\psi$ compares with the effectiveness of Rasmussen's $s$-invariant~\cite{Ras10} and the self-linking number $\self$ combined. Before we prove \fullref{thm:effectiveness}, we recall what is known about the relationship among $s$, $\self$, and $\psi$. For any transverse knot $K$, Plamenevskaya shows that $\self(K)\leq s(K)-1$ \cite[Proposition~4]{olga-psi}. Further, Baldwin and Plamenevskaya show that if this is an equality, then $\psi(K)\neq 0$ \cite[Theorem~1.2]{BP-Kh}. They also show that the converse holds for some knots, but the question of whether there is a knot $K$ with $s(K) > \self(K)+1$ and $\psi(K)\not=0$ is still open.

Let $\Leg_-$ and $\Leg_+$ be two Legendrian knots with positive transverse pushoffs $K_-$ and $K_+$. Suppose $\psi(K_{\pm})\neq 0$ if and only if $s(K_{\pm}) = \self(K_{\pm})+1$. Then \fullref{cor:Vanish} will not provide a stronger obstruction to the existence of a decomposable Lagrangian cobordism from $\Leg_-$ to $\Leg_+$ than a combination of classical Legendrian invariants and the $s$-invariant.

\begin{proof}[Proof of \fullref{thm:effectiveness}]

    Suppose that we wish to obstruct  the existence of a  decomposable Lagrangian cobordism with Euler characteristic $E$ from $\Leg_-$ to $\Leg_+$.   Suppose further that the classical Legendrian invariants $\tb$ and $r$ alone  do not obstruct the existence of such a cobordism. In this case, we have 
    \begin{align*}
    \rot(\Leg_-) &= \rot(\Leg_+),\\
    \tb(\Leg_-) &= \tb(\Leg_+) + E.
    \end{align*}

    Consider the positive transverse pushoffs $K_-$ and $K_+$ of  $\Leg_-$ and $\Leg_+$. We have
    \[
    \self(K_+) =\tb(\Leg_+) - \rot(\Leg_+) = \tb(\Leg_-)- E - \rot(\Leg_-)  = \self(K_-) - E.
    \]

    Now suppose $\psi(K_-) \not=0$ and $\psi(K_+) = 0$.
    The fact that $\psi(K_+)=0$ guarantees that $\self(K_+) < s(K_+)-1$, so 
    \begin{equation}\label{eqn:k1}
    \self(K_-) -  E < s(K_+)-1.
    \end{equation}
    Now, if $K_-$ is a transverse knot where $\psi(K_-) \neq 0$ if and only if $\self(K_-) = s(K_-)-1$, then 
     $\psi(K_-)\not=0$ implies 
    \begin{equation}\label{eqn:k2}
    \self(K_-) = s(K_-)-1,
    \end{equation}
    and  \fullref{eqn:k1} and \fullref{eqn:k2} combined imply that 
    \[
    s(K_-)-E < s(K_+).
    \]
    This difference between $s(K_-)$ and $s(K_+)$ obstructs the existence of a connected smooth cobordism with Euler characteristic $E$ from $K_-$ to $K_+$. Since decomposable Lagrangian cobordisms between knots are automatically connected, the Rasmussen $s$ invariant, in combination with $\tb$ and $r$,  already obstructs the existence of a decomposable Lagrangian cobordism  from $\Leg_-$ to $\Leg_+$.
\end{proof}

The third ingredient of Ekholm, Honda, and K\'{a}lm\'{a}n's question involves what we refer to as ``Ng's line'' of a knot $K$. Before we prove \fullref{thm:line}, we recall the definition of this line and reconcile notational conventions from \cite{NgLenhard2005ALTb} and \cite{EHK16}.

In \cite{NgLenhard2005ALTb}, Ng observed an interesting relationship between Khovanov homology and the Thurston-Bennequin number. Specifically, if $\overline{\tb}(K)$ denotes the maximal Thurston-Bennequin number over all Legendrian links in the smooth link type of $K$, then Ng proved that 
    \[\overline{\tb}(K) \leq \min \{ k  \mid \bigoplus_{i-j=k}\Kh^{i,j}(K)\neq 0\};\]
here, $i$ is the quantum grading and $j$ is the homological grading. Denote the right hand side of Ng's inequality by $C$. Geometrically, one can think of the equation $i-j=C$ as the highest line of slope one on the $ji$-plane such that Khovanov homology is supported on and above that line; this is what became colloquially known as ``Ng's line for $K$.''

In \cite[Example 9.6]{EHK16}, Ekholm, Honda, and K\'{a}lm\'{a}n  observed that when $\Leg$ is the right-handed trefoil, the elements of $\LKh \coloneqq \{F_{m L}(\pm 1)\}\subset \Kh(m(\Leg))$ are supported on the dual of Ng's line for $\Leg$, under the mirroring duality $\Kh^{i,j}(\Leg)\otimes_{\Z}\Q \cong \Kh^{-i,-j}(m(\Leg))\otimes_{\Z}\Q$ from \cite[Corollary~11]{Kh00}. We caution the reader that this isomorphism is used implicitly throughout \cite{EHK16}, and that the grading conventions in \cite{EHK16} are opposite to those in \cite{NgLenhard2005ALTb}. Specifically, when Ekholm, Honda, and K\'{a}lm\'{a}n discuss the line $j-i = C$ where $C = \max \{ k  \mid \bigoplus_{j-i=k}\Kh^{i,j}(m(K))\neq 0\}$ in \cite{EHK16}, they use $i$ to denote the homological grading and $j$ to denote the quantum grading. Further, when they refer to this line as ``Ng's line'', they mean neither Ng's line for $m(K)$, nor Ng's line for $K$, but rather the image of Ng's line for $K$ under the plane rotation $(x,y)\mapsto (-x,-y)$.

Returning to the classical invariants, we can see that the observation in \cite[Example 9.6]{EHK16} is not a coincidence, but always the case for knots with decomposable Lagrangian fillings.

\begin{proof}[Proof of \fullref{thm:line}]
    Let $L$ be a decomposable Lagrangian filling for $\Leg$, and let $K$ denote the positive transverse pushoff of $\Leg$. 

    Since $\Leg$ admits an exact Lagrangian filling, we know that $r(\Leg)=0$, and so $\self(K) = \tb(\Leg) = -\chi(L)$. Thus, the bigrading of $\psi(K)$ is $(-\chi(L),0)$. From \fullref{thm:KhLagrangian}, we know that $\psi(K)\neq 0$, so $\psi(K)$ lies on or above Ng's line, i.e.\ $\tb(\Leg) = -\chi(L) \geq \min \{ i-j  \mid \Kh^{i,j}(K)\neq 0\}$. Combining with Ng's inequality, we see that $tb(\Leg)\geq \overline{\tb}(K)$, so in fact $tb(\Leg) = \overline{\tb}(K) = -\chi(L)$ and thus Ng's line for $K$ is the line $i-j = \overline{\tb}(K)$, and $\psi(K)$ lies on this line.

    Now, the map $F_{mL}$ from \cite{J04} is homogeneous of bidegree $(\chi(L), 0)$, and we saw that $\chi(L) = -\tb(\Leg)$, so $\LKh$ is a set of non-torsion homogeneous elements in $\Kh^{-\tb(\Leg), 0}(m(K))$. Thus, the images of these elements under the mirroring duality lie on Ng's line for $K$.
\end{proof}

\fullref{cor:top-obstr} follows immediately, using the above observation that if $\psi$ is nonzero for the positive transverse pushoff of $\Leg$, then Ng's line for $K$ is the line $\overline{\tb}(K) = i-j$.

\section{Filtered invariants and decomposable Lagrangian cobordisms}

In~\cite{LipshitzRobert2015Otif}, Lipshitz, Ng, and Sarkar introduce a generalization of Plamenevskaya's invariant $\psi(K)$ utilizing Bar-Natan's deformation of Khovanov homology. They note that their construction also works in Lee's deformation $\Lee(K)$ \cite{Lee02}, which we will summarize below.

Let $\CLee(K)$ denote the Lee complex associated to a link diagram $K$, equipped with a homological grading and a filtration by the quantum grading. To each orientation $o$ of $K$, we can associate a generator $\psigen(o) \in \CLee(K)$. Similar to the generator of $\psi(K)$, $\psigen(o)$ lives in the resolution in which every crossing is replaced by an oriented smoothing. To each circle $C$ in this oriented smoothing, we associate the sum or difference of the generators $v_+$ and $v_-$ corresponding to that circle. More precisely, given a circle $C$, pick a point $c$ on $C$. Push $c$ to the left of $C$ (relative to the orientation) to get a point $c'$, then choose an arc from $c'$ to the point at infinity in the exterior of the diagram. If this arc crosses the circles in the oriented smoothing an even number of times, then we associate to $C$ the sum $v_- + v_+$. Otherwise, we associate to $C$ the difference $v_- - v_+$. The generator $\psigen(o)$ is then the tensor product of the factors associated to each circle $C$.

The filtered Plamenevskaya invariant $\psigen(K)$ of a transverse link $K$ is defined to be the generator $\psigen(o)$ corresponding to the positive orientation $o$ of $K$, up to filtered homotopy equivalence. Lipshitz, Ng, and Sarkar further define invariants $\psigen_{p,q}(K)$ for $p \le 0 < q$, which live in various subquotients of the Lee complex. For brevity, given a chain complex $\complex$, let $\gr_{p,q}(\complex)$ denote the subquotient complex
\begin{equation*}
    \gr_{p,q}(\complex) = \mathcal{F}_{sl(K)+2p} \, \complex / \mathcal{F}_{sl(K)+2q} \, \complex \,.
\end{equation*}
The filtered invariants $\psigen_{p,q}(K)$ live in the complex
\begin{equation*}
    \psigen_{p,q}(K) \in \gr_{p,q}(\CLee^0(K)) \,.
\end{equation*}
These invariants are defined as the images of $\psigen(K)$ under the implied quotients. These filtered invariants also recover the original Plamenevskaya invariant. Since the associated graded object of the Lee complex $\CLee(K)$ is the Khovanov complex $\CKh(K)$, the image of $\psigen(o)$ in the associated graded complex is the generator of the Plamenevskaya invariant $\psi(K)$.

In particular, $\psigen_{0,1}(K)$ generates $\psi(K)$. The invariants $\psigen_{p,q}(K)$ also interact nicely with the cobordism maps from ascending cobordisms with positive critical points; this is analogous to the situation with the spectral invariants from grid homology~\cite{JPSWW22}.

\begin{theorem}\label{thm:Spectral}
     Let $\Leg_-$ and $\Leg_+$ be two Legendrian links with positive transverse pushoffs $K_-$ and $K_+$, and suppose there is a decomposable Lagrangian cobordism $L$ from $\Leg_-$ to $\Leg_+$. Then the induced cobordism maps \[ F_L^{p,q}\colon \gr_{p,q}(\CLee^0(\Leg_+)) \to \gr_{p,q}(\CLee^0(\Leg_-)) \]
     send $[\psigen_{p,q}(K_+)]$ to $\pm[\psigen_{p,q}(K_-)]$.
\end{theorem}

Just as in the non-filtered case we have the following immediate corollary to \fullref{thm:Spectral}.

\begin{corollary}
\label{cor:SpectralVanish}
    Let $\Leg_-$ and $\Leg_+$ be two Legendrian links with positive transverse pushoffs $K_-$ and $K_+$. If for some $0\leq p <q$ we have $[\psigen_{p,q}(K_-)] \not=0$ and $[\psigen_{p,q}(K_+)] = 0$, then there is no decomposable Lagrangian cobordism $L$ from $\Leg_-$ to $\Leg_+$.
\end{corollary}

Unlike \fullref{cor:Vanish}, there is a possibility that \fullref{cor:SpectralVanish} provides an interesting obstruction even if $\psi$ is always determined by the $s$-invariant and the self-linking number. However, there are not many known computations of $\psigen_{p,q}$ at the moment, and no computer program to compute it to the best of our knowledge, so we can not say for sure if \fullref{cor:SpectralVanish} can provide new obstructions to decomposable Lagrangian cobordisms.

The proof of \fullref{thm:Spectral} relies on algebraic lemmas about maps between filtered complexes from~\cite{JPSWW22} and an brief analysis of maps induced by the elementary cobordisms used to build ascending cobordisms with positive critical points.

\begin{proof}[Proof of \fullref{thm:Spectral}]
    We will show that on the chain level the cobordism map associated to an ascending cobordism with positive critical points satisfies the hypotheses of~\cite[Lemma 2.4]{JPSWW22}. With the work from~\cite{LipshitzRobert2015Otif} showing the classes are transverse invariants, all that is left to consider are birth and the saddle cobordisms.

    If the cobordism $L$ is a birth, then, considering the map from $K_+$ to $K_-$, we see a death of a trivial circle. In the definition of $\psigen_{p,q}$, every circle in the oriented resolution is labeled by $v_- + v_+$ or $v_- - v_+$, both of which are sent to $1$ under the map induced by the cap cobordism. Therefore, $\psigen_{p,q}(K_+)$ is mapped to $\psigen_{p,q}(K_-)$ on the chain level.

    If $L$ is a saddle, then the map induced on Khovanov homology can be realized as an isotopy followed by a positive crossing resolution, as depicted (in reverse) in \fullref{Fig:Band}. A generalization of Plamenevskaya's argument~\cite{olga-psi} for the behavior of her invariant under the resolution of a positive crossing to a braid applies. Plamenevskaya argues that if both strands of the oriented resolution are labeled with $v_-$, then under the cobordism map consisting of a split followed by a cap, the image also contains both strands labeled with $v_-$. It is true more generally that any labeling of the two strands of the oriented resolution with $v_+$ or $v_-$ will be preserved by this map. Since each strand in $\psigen(o)$ is labeled by $v_- + v_+$ or $v_- - v_+$, linearity ensures that $\psigen(o)$ is also preserved by this cobordism. Therefore, $\psigen_{p,q}(K_+)$ is mapped to $\psigen_{p,q}(K_-)$ on the chain level.
\end{proof}

\bibliographystyle{mwamsalphack}
\bibliography{references}

@article {BLW22,
    AUTHOR = {Baldwin, John A. and Lidman, Tye and Wong, C.-M. Michael},
     TITLE = {Lagrangian cobordisms and {L}egendrian invariants in knot
              {F}loer homology},
   JOURNAL = {Michigan Math. J.},
  FJOURNAL = {Michigan Mathematical Journal},
    VOLUME = {71},
      YEAR = {2022},
    NUMBER = {1},
     PAGES = {145--175},
      ISSN = {0026-2285},
   MRCLASS = {57K10 (57K33)},
  MRNUMBER = {4389674},
       DOI = {10.1307/mmj/20195786},
       URL = {https://doi.org/10.1307/mmj/20195786},
}

@article {BS18,
    AUTHOR = {Baldwin, John A. and Sivek, Steven},
     TITLE = {Invariants of {L}egendrian and transverse knots in monopole
              knot homology},
   JOURNAL = {J. Symplectic Geom.},
  FJOURNAL = {The Journal of Symplectic Geometry},
    VOLUME = {16},
      YEAR = {2018},
    NUMBER = {4},
     PAGES = {959--1000},
      ISSN = {1527-5256},
   MRCLASS = {57M27 (57R58)},
  MRNUMBER = {3917725},
MRREVIEWER = {Daniel Ruberman},
       DOI = {10.4310/JSG.2018.v16.n4.a3},
       URL = {https://doi.org/10.4310/JSG.2018.v16.n4.a3},
}

@article {BS21,
    AUTHOR = {Baldwin, John A. and Sivek, Steven},
     TITLE = {On the equivalence of contact invariants in sutured {F}loer
              homology theories},
   JOURNAL = {Geom. Topol.},
  FJOURNAL = {Geometry \& Topology},
    VOLUME = {25},
      YEAR = {2021},
    NUMBER = {3},
     PAGES = {1087--1164},
      ISSN = {1465-3060},
   MRCLASS = {53D40 (53D10 57R58)},
  MRNUMBER = {4268162},
MRREVIEWER = {Jun Zhang},
       DOI = {10.2140/gt.2021.25.1087},
       URL = {https://doi.org/10.2140/gt.2021.25.1087},
}

@inproceedings {CDGG15,
    AUTHOR = {Chantraine, Baptiste and Dimitroglou Rizell, Georgios and
              Ghiggini, Paolo and Golovko, Roman},
     TITLE = {Floer homology and {L}agrangian concordance},
 BOOKTITLE = {Proceedings of the {G}\"{o}kova {G}eometry-{T}opology {C}onference
              2014},
     PAGES = {76--113},
 PUBLISHER = {G\"{o}kova Geometry/Topology Conference (GGT), G\"{o}kova},
      YEAR = {2015},
   MRCLASS = {53D42 (53D12 57R17 57R58)},
  MRNUMBER = {3381440},
MRREVIEWER = {Daniel V. Mathews},
Customkey = {CDGG},
}

@article {Cha15:NotSym,
    AUTHOR = {Chantraine, Baptiste},
     TITLE = {Lagrangian concordance is not a symmetric relation},
   JOURNAL = {Quantum Topol.},
  FJOURNAL = {Quantum Topology},
    VOLUME = {6},
      YEAR = {2015},
    NUMBER = {3},
     PAGES = {451--474},
      ISSN = {1663-487X},
   MRCLASS = {57R17 (53D42 57M27 57M50)},
  MRNUMBER = {3392961},
MRREVIEWER = {Daniel V. Mathews},
       DOI = {10.4171/QT/68},
       URL = {https://doi.org/10.4171/QT/68},
}

@article {CNS16,
    AUTHOR = {Cornwell, Christopher and Ng, Lenhard and Sivek, Steven},
     TITLE = {Obstructions to {L}agrangian concordance},
   JOURNAL = {Algebr. Geom. Topol.},
  FJOURNAL = {Algebraic \& Geometric Topology},
    VOLUME = {16},
      YEAR = {2016},
    NUMBER = {2},
     PAGES = {797--824},
      ISSN = {1472-2747},
   MRCLASS = {57M25 (53D12 53D42 57R17)},
  MRNUMBER = {3493408},
MRREVIEWER = {Georgios Dimitroglou Rizell},
       DOI = {10.2140/agt.2016.16.797},
       URL = {https://doi.org/10.2140/agt.2016.16.797},
}

@article {EHK16,
    AUTHOR = {Ekholm, Tobias and Honda, Ko and K\'{a}lm\'{a}n, Tam\'{a}s},
     TITLE = {Legendrian knots and exact {L}agrangian cobordisms},
   JOURNAL = {J. Eur. Math. Soc. (JEMS)},
  FJOURNAL = {Journal of the European Mathematical Society (JEMS)},
    VOLUME = {18},
      YEAR = {2016},
    NUMBER = {11},
     PAGES = {2627--2689},
      ISSN = {1435-9855},
   MRCLASS = {53D42 (53D10 53D40 57M27 57R90)},
  MRNUMBER = {3562353},
MRREVIEWER = {Georgios Dimitroglou Rizell},
       DOI = {10.4171/JEMS/650},
       URL = {https://doi.org/10.4171/JEMS/650},
}

@article {GJ19,
    AUTHOR = {Golla, Marco and Juh\'{a}sz, Andr\'{a}s},
     TITLE = {Functoriality of the {EH} class and the {LOSS} invariant under
              {L}agrangian concordances},
   JOURNAL = {Algebr. Geom. Topol.},
  FJOURNAL = {Algebraic \& Geometric Topology},
    VOLUME = {19},
      YEAR = {2019},
    NUMBER = {7},
     PAGES = {3683--3699},
      ISSN = {1472-2747},
   MRCLASS = {57K18 (57K33)},
  MRNUMBER = {4045364},
       DOI = {10.2140/agt.2019.19.3683},
       URL = {https://doi.org/10.2140/agt.2019.19.3683},
}

@article {JPSWW22,
    Author = {Jubeir, Mitchell and Petkova, Ina and Schwartz, Noah and Winkeler, Zachary and Wong, C.-M. Michael},
    Title = {Spectral {GRID} invariants and {L}agrangian cobordisms},
    journal = "Studia Scientiarum Mathematicarum Hungarica",
      year = "2025",
      publisher = "Akad{\'e}miai Kiad{\'o}",
      address = "Budapest, Hungary",
      volume = "61",
      number = "4",
      doi = "10.1556/012.2024.04325",
      pages=      "311 - 340",
      url = "https://akjournals.com/view/journals/012/61/4/article-p311.xml"
}

@article{Kh00,
author = {Mikhail Khovanov},
title = {{A categorification of the Jones polynomial}},
volume = {101},
journal = {Duke Mathematical Journal},
number = {3},
publisher = {Duke University Press},
pages = {359 -- 426},
year = {2000},
doi = {10.1215/S0012-7094-00-10131-7},
URL = {https://doi.org/10.1215/S0012-7094-00-10131-7}
}

@article {Pan17,
    AUTHOR = {Pan, Yu},
     TITLE = {The augmentation category map induced by exact {L}agrangian
              cobordisms},
   JOURNAL = {Algebr. Geom. Topol.},
  FJOURNAL = {Algebraic \& Geometric Topology},
    VOLUME = {17},
      YEAR = {2017},
    NUMBER = {3},
     PAGES = {1813--1870},
      ISSN = {1472-2747},
   MRCLASS = {53D42 (53D12 57M50 57R17)},
  MRNUMBER = {3677941},
MRREVIEWER = {Paolo Ghiggini},
       DOI = {10.2140/agt.2017.17.1813},
       URL = {https://doi.org/10.2140/agt.2017.17.1813},
}

@article {ST13,
    AUTHOR = {Sabloff, Joshua M. and Traynor, Lisa},
     TITLE = {Obstructions to {L}agrangian cobordisms between {L}egendrians
              via generating families},
   JOURNAL = {Algebr. Geom. Topol.},
  FJOURNAL = {Algebraic \& Geometric Topology},
    VOLUME = {13},
      YEAR = {2013},
    NUMBER = {5},
     PAGES = {2733--2797},
      ISSN = {1472-2747},
   MRCLASS = {53D12 (57Q60 57R17)},
  MRNUMBER = {3116302},
MRREVIEWER = {Jelena Kati\'{c}},
       DOI = {10.2140/agt.2013.13.2733},
       URL = {https://doi.org/10.2140/agt.2013.13.2733},
}

@article{hayden_khovanov_2024,
	title = {Khovanov homology and exotic surfaces in the 4-ball},
	volume = {2024},
	copyright = {De Gruyter expressly reserves the right to use all content for commercial text and data mining within the meaning of Section 44b of the German Copyright Act.},
	issn = {1435-5345},
	url = {https://www.degruyter.com/document/doi/10.1515/crelle-2024-0001/pdf},
	doi = {10.1515/crelle-2024-0001},
	language = {en},
	number = {809},
	urldate = {2024-07-02},
	journal = {Journal für die reine und angewandte Mathematik (Crelles Journal)},
	author = {Hayden, Kyle and Sundberg, Isaac},
	month = apr,
	year = {2024},
	pages = {217--246},
}

@article{olga-psi,
author = {Plamenevskaya, Olga},
year = {2005},
month = {01},
pages = {571--586},
title = {Transverse knots and {K}hovanov homology},
volume = {13},
journal = {Mathematical Research Letters},
doi = {10.4310/MRL.2006.v13.n4.a7}
}

@article{BP-Kh,
author = {John A. Baldwin and Olga Plamenevskaya},
title = {Khovanov homology, open books, and tight contact structures},
journal = {Advances in Mathematics},
volume = {224},
number = {6},
pages = {2544-2582},
year = {2010},
issn = {0001-8708},
doi = {https://doi.org/10.1016/j.aim.2010.02.010},
url = {https://www.sciencedirect.com/science/article/pii/S0001870810000654}
}

@article{NgLenhard2005ALTb,
issn = {1472-2739},
journal = {Algebraic \& geometric topology},
pages = {1637--1653},
volume = {5},
publisher = {Geometry & Topology Publications},
number = {4},
year = {2005},
title = {A {L}egendrian {T}hurston-{B}ennequin bound from {K}hovanov homology},
language = {eng},
address = {COVENTRY},
author = {Ng, Lenhard},
}

@article{LipshitzRobert2015Otif,
issn = {1663-487X},
journal = {Quantum topology},
pages = {475--513},
volume = {6},
publisher = {European Mathematical Society Publishing House},
number = {3},
year = {2015},
title = {On transverse invariants from {K}hovanov homology},
copyright = {European Mathematical Society},
language = {eng},
address = {Zuerich, Switzerland},
author = {Lipshitz, Robert and Ng, Lenhard and Sarkar, Sucharit},
}

@article{J04,
   title={An invariant of link cobordisms from {K}hovanov homology},
   volume={4},
   ISSN={1472-2747},
   url={http://dx.doi.org/10.2140/agt.2004.4.1211},
   DOI={10.2140/agt.2004.4.1211},
   number={2},
   journal={Algebraic \& Geometric Topology},
   publisher={Mathematical Sciences Publishers},
   author={Jacobsson, Magnus},
   year={2004},
   month=dec, pages={1211–1251} }

@article{Ras10,
author = {Rasmussen, Jacob},
year = {2004},
month = {03},
pages = {419--447},
title = {Khovanov homology and the slice genus},
volume = {182},
journal = {Inventiones mathematicae},
doi = {10.1007/s00222-010-0275-6}
}

@article{SundbergIsaac2022RKc,
issn = {1472-2739},
journal = {Algebraic \& geometric topology},
pages = {3983--4008},
volume = {22},
publisher = {Geometry & Topology Publications},
number = {8},
year = {2022},
title = {Relative {K}hovanov-{J}acobsson classes},
language = {eng},
address = {COVENTRY},
author = {Sundberg, Isaac and Swann, Jonah},
}

@ARTICLE{Lee02,
       author = {{Lee}, Eun Soo},
        title = "{An endomorphism of the Khovanov invariant}",
      journal = {arXiv Mathematics e-prints},
         year = 2002,
        month = oct,
          eid = {math/0210213},
        pages = {math/0210213},
          doi = {10.48550/arXiv.math/0210213},
archivePrefix = {arXiv},
       eprint = {math/0210213},
 primaryClass = {math.GT},
       adsurl = {https://ui.adsabs.harvard.edu/abs/2002math.....10213L}
}

\end{document}